\documentstyle[12pt]{article}
\begin{document}
\baselineskip=15pt
\begin{center}
{\Large{\bf Analogy between arithmetic of elliptic curves and conics}}\\
\vspace{0.5cm}
H. Gopalkrishna Gadiyar and R. Padma\\
AU-KBC Research Centre, M. I. T. Campus, Anna University\\
Chromepet, Chennai 600 044 INDIA\\
E-mail:gadiyar@au-kbc.org, padma@au-kbc.org
\end{center}

\vspace{1cm}

\begin{center}
{\bf Abstract}
\end{center}

In this brief note we bring out the analogy between the arithmetic of elliptic curves and the Riemann zeta-function. 
\vspace{1cm}
$$ 
\cdots 
$$

The aim of this note is to point out the possibility of developing the theory of elliptic curves so that the arithmetical and analytic aspects are developed in strict analogy with the classical theory of the Riemann zeta-function. This has been triggered by Lemmermeyer's perceptive analysis of conics [3-6]. 

The Riemann zeta-function has two fundamental representations. 
\begin{equation}
\sum_{n=1}^\infty \frac{1}{n^s} ~=~ \prod_p \left(1-\frac{1}{p^s}\right)^{-1}, ~if~Re. ~s~>~1 \,.
\end{equation}
These we call the left hand and the right hand side. The left hand side is a sum over the integers and the right hand side is a product over primes.The left hand side is got from trigonometry. 
It begins with Euler's expansion of $\sin x$ into an infinite product given by the formula
\begin{equation}
\sin x ~=~ x \prod_n (1 - \frac{x^2}{n^2\pi ^2}) \, ,
\end{equation}
and the expression in terms of the exponentials\footnote{Euler got the value of $\zeta (2)$ by expanding the sine function in Taylor series and comparing coefficients of $x^3$ on both sides}: 
\begin{equation}
\sin x ~=~ \frac{e^{ix} - e^{-ix}}{2i} \, . 
\end{equation}
Next Euler took the logarithmic derivative and got an equation where one side consisted of fractions 
\begin{equation}
\frac{1}{x} + 2 \sum_{n=1}^\infty \frac{x}{x^2 - n^2\pi ^2}
\end{equation}
while the other side depended on exponentials. 
Riemann essentially interpolated this equation by integrating both sides from $0$ to $\infty$ after multiplying by $x^{s-1}$ (in other words, Mellin transform) and obtained the functional equation for the Riemann zeta-function $\zeta (s)$(See Appendix 2.)

In the case of elliptic curves [1,2,8,10], the theory of the subject falls into two parts. A classical part consisting of developing the theory in strict analogy with trigonometric functions ($\wp(z) , \zeta (z)$ and $\sigma (z)$ are the analogs of $\csc ^2z, ~\cot z,$~and $\sin z$. Elliptic functions reduce to these trigonometric functions in the singular limit.) Then the modern number theoretic aspects have been developed where elliptic curves are studied modulo $p$, $L$ - functions constructed and pasted together yielding 
\begin{equation}
L(s,E)~=~\prod_{p|\Delta} \left(\frac{1}{1-a_pp^{-s}}\right) \prod_{p \mid\!\!\!/ \Delta}\left(\frac{1}{1-a_p p^{-s}+p^{1-2s}}\right) 
\end{equation}
The exercises in P\'{o}lya  and Szeg\"{o} (see Appendix 1) seem to suggest that for extracting the arithmetic an ansatz of the form 
\begin{equation}
\prod _{n=1}^{\infty }( 1-\frac{x^2}{n^2 \pi ^2})^{A_n}
\end{equation}
would exist for $\sigma (z)$ in analogy with 
\begin{equation}
\sin x ~=~ x \prod_n (1 - \frac{x^2}{n^2\pi ^2}) \, ,
\end{equation}
This would lead to corresponding formulae for $\zeta (z)$ and $\wp (z)$. 

Recall that the Shimura-Taniyama-Weil conjecture (see Appendix 3.) is written in the language of modular functions. The entire conjecture can be recast into the trigonometric form. 
All essential arithmetic information is encoded in the elliptic functions. An ansatz of the form (in a formal power series sense)
\begin{equation}
\zeta (z) = \frac{1}{z} + \sum_n \frac{2 a(n) z}{n^2+z^2}
\end{equation}
gives a solution to the equation
\begin{equation}
4 (\wp ^\prime (z))^2 = 4 \wp ^3(z) - g_2 \wp (z) -g_3
\end{equation}
Corresponding to the map $\tau \rightarrow -1/\tau $ (see Appendix 3.) is the classical mapping (Fourier transform) used in the Poisson summation formula of going from partial fractions to exponentials. This can be easily seen in the case of $\sin x$. $N$ is a scale for $z$. The exercises given in Appendix 1. show that the trigonometric form and the modular form are related via Lambert series by the familiar trick of logarithmic derivative. 

The efforts of the last two decades have been spent in successfully exploiting the right hand side. This has led to significant advances [11a,11b]. The analysis carried out in this note shows that the left hand side needs to be developed. That is, a completely analytic approach to elliptic curves needs to be developed. The classical theory of elliptic functions has to be reworked to extract the arithmetical properties. 

\begin{center}
{\bf Appendix 1.}
\end{center}

\noindent See Exercises 65,66,67,68 in [7].\\
{\bf 65.} From the identity
\begin{equation}
\zeta (s) \sum_{n=1}^\infty a_n n^{-s} = \sum_{n=1}^\infty A_n n^{-s}
\end{equation}
it follows that
\begin{equation}
\sum_{n=1}^\infty \frac{a_n x^n}{1-x^n}=\sum_{n=1}^\infty A_n x^n
\end{equation}
and conversely, where $a_1, a_2, a_3,...,A_1, A_2, A_3,...$ are constants. (The series occurring on the left hand side of the second equation is called a {\it Lambert series}.

\noindent Proof: From the first identity we obtain
\begin{equation}
A_n= \sum_{t|n} a_t \, ,
\end{equation}
and from the second identity
\begin{equation}
\sum_{n=1}^\infty \frac{a_n x^n}{1-x^n}= \sum_{m=1}^\infty a_m (x^m+x^{2m}+x^{3m}+....)=\sum_{n=1}^\infty (\sum_{t|n}a_t)x^n \, .
\end{equation}

\noindent {\bf 66.} The two identities
\begin{equation}
\zeta (s) (1-2^{1-s})\sum_{n=1}^\infty a_n n^{-s} = \sum_{n=1}^\infty B_n n^{-s} \, ,
\end{equation}
\begin{equation}
\sum_{n=1}^\infty \frac{a_n x^n}{1+x^n}=\sum_{n=1}^{\infty} B_n x^n
\end{equation}
are equivalent.

\noindent Proof: We have
\begin{equation}
(1-2^{1-s})\zeta (s) = 1^{-s} - 2^{-s} +3^{-s}-4^{-s}+...
\end{equation}
From the first identity it therefore follows that
\begin{equation}
B_n =  \sum_{t|n} (-1)^{t-1} a_{\frac{n}{t}} \, ,
\end{equation}
and from the second identity 
\begin{equation}
\sum_{n=1}^\infty \frac{a_n x^n}{1+x^n}= \sum_{m=1}^\infty a_m (x^m-x^{2m}+x^{3m}-....)=\sum_{n=1}^\infty (\sum_{t|n}(-1)^{t-1} a_{\frac{n}{t}})x^n \, .
\end{equation}

\noindent {\bf 67.} Assume that between the numbers $a_1, a_2, a_3,...$ and $A_1, A_2, A_3,...$ the same relation subsists as in {\bf 65}. Then
\begin{equation}
\prod_{n=1}^\infty ( \frac{n}{x} \sin \frac{x}{n} )^{a_n} = \prod _{n=1}^{\infty }( 1-\frac{x^2}{n^2 \pi ^2})^{A_n}
\end{equation}
Proof: Writing $A_n$ in the form $A_n=\sum_{t|n} a_t$ and collecting the factors with exponent $a_m$, we obtain on the right hand side the product
\begin{equation}
(1-\frac{x^2}{m^2\pi ^2})^{a_m} (1-\frac{x^2}{(2m)^2\pi ^2})^{a_m} (1-\frac{x^2}{(3m)^2\pi ^2})^{a_m}... \, .
\end{equation}
On the left hand side we have correspondingly the $a_m$th power of
\begin{equation}
\frac{m}{x} \sin \frac{x}{m}=(1-\frac{x^2}{m^2\pi ^2}) (1-\frac{x^2}{(2m)^2\pi ^2})(1-\frac{x^2}{(3m)^2\pi ^2})... \, .
\end{equation}

\noindent {\bf 68.} Assume that between the numbers $a_1, a_2, a_3,...$ and $B_1, B_2, B_3,...$ the same relation subsists as in {\bf 66}.  Then
\begin{equation}
\prod_{n=1}^\infty ( \frac{x}{2n} \cot \frac{x}{2n} )^{a_n} = \prod _{n=1}^{\infty }( 1-\frac{x^2}{n^2 \pi ^2})^{B_n} \, .
\end{equation}
Proof:
\begin{equation}
\prod_{k=1}^\infty \left(1-\frac{x^2}{(km)^2\pi ^2}\right)^{(-1)^{k-1}}=\frac{x}{2m} \cot \frac{x}{2m} \, .
\end{equation}

\begin{center}
{\bf Appendix 2.}
\end{center}

This method [9] uses the partial fraction expansion of $\coth $ function which is given by:
\begin{equation}
\frac{1}{e^x -1} = \frac{1}{x} - \frac{1}{2} + 2x \sum_{n=1}^\infty \frac{1}{x^2 + 4n^2\pi ^2}
\end{equation}
which when applied to the integral
\begin{equation}
\Gamma (s) \zeta (s)~=~ \int_0^\infty \left( \frac{1}{e^x -1} - \frac{1}{x} + \frac{1}{2}\right) x^{s-1} dx \, , (-1~<~Re.~ s~<0) \, ,
\end{equation}
gives the functional equation of $\zeta (s)$.
\begin{eqnarray}
\Gamma (s) \zeta (s)&=& \int_0^\infty 2x \sum_{n=1}^\infty \frac{1}{x^2 +4 n^2\pi ^2} ~x^{s-1} dx \\
~&=& 2 \sum_{n=1}^\infty \int_0^\infty \frac{x^s}{x^2 +4 n^2\pi ^2} dx \\
~&=& 2 \sum_{n=1}^\infty (2n\pi)^{s-1} \frac{\pi}{2 \cos \frac{1}{2}s\pi}\\
~&=& \frac{2^{s-1}\pi ^s}{\cos \frac{1}{2}s\pi} \zeta (1-s) \,.
\end{eqnarray}

\begin{center}
{\bf Appendix 3.}
\end{center}

Let $E/\cal{Q}$ be an elliptic curve of conductor $N$, let $L_E(s)~=~\sum_{n} c_n n^{-s}$ be its $L$-series, and let $f(\tau)~=~\sum c_n~e^{2\pi in\tau} $ be the inverse Mellin transform of $L_E$.\\
(a) $f(\tau)$ is a cusp form of weight 2 for the congruence subgroup $\Gamma _0(N)$ of $SL_2(Z)$.\\
(b) For each prime $p \mid\!\!\!\!/ N$, let $T(p)$ be the corresponding Hecke operator; and let $W$ be the operator $(Wf)(\tau )~=~f(-1/N\tau)$. Then
\begin{equation}
T(p)f~=~c_p~f ~~\mbox{and}~Wf~=~wf \, ,
\end{equation}
where $w~=~\pm 1$ is the sign of the functional equation
\begin{equation}
\Lambda (s,E)~=~ w~\Lambda (2-s,E)
\end{equation}
where 
\begin{equation}
\Lambda (s,E)~=~N^{\frac{s}{2}}~(2\pi)^{-s} \Gamma (s)~L_E(s) \, .
\end{equation}
(c) Let $\omega $ be an invariant differential on $E/\cal{Q}$. There exists a morpism $\phi : X_0(N) \rightarrow E$, defined over $\cal{Q}$, such that $\phi ^*(\omega)$ is a multiple of the differential form on $X_0(N)$ represented by $f(\tau )~d\tau $.

\begin{center}
{\bf References}
\end{center}

\begin{description}
\item[1.] A. W. Knapp, Elliptic Curves, Princeton University Press, 1993.
\item[2.] N. Koblitz, Introduction to Ellitpic Curves and Modular Forms, GTM 97, Springer - Verlag, 1984.
\item[3.] F. Lemmermeyer, Conics - A poor man's elliptic curves, \\arXiv:Math.NT/0311306v1, 18 Nov 2003.
\item[4.] F. Lemmermeyer, Higher descent of Pell conics. I. From Legendre to Selmer, arXiv:Math.NT/0311309v1, 18 Nov 2003.
\item[5.] F. Lemmermeyer, Higher descent of Pell conics. II. Two centuries of missed opportunities, arXiv:Math.NT/0311296v1, 18 Nov 2003.
\item[6.] F. Lemmermeyer, Higher descent of Pell conics. III. The first 2 - descent, arXiv:Math.NT/0311310v1, 18 Nov 2003.
\item[7.] G. P\'{o}lya  and G. Szeg\"{o}, Problems and Theorems in Analysis, Volume II, Springer International Student Edition, 1976.
\item[8.] J.H. Silverman, The Arithmetic of Elliptic Curves, GTM 106, Springer - Verlag, 1985.
\item[9.] E. C. Titchmarsh, Revised by D. R. Heath - Brown, The theory of the Riemann Zeta-function, Second Edition, Oxford Science Publications, 1986.
\item[10.] Andre Weil, Elliptic Functions according to Eisenstein and Kronecker, Ergebnisse der Mathematik und ihrer Grenzgebiete 88, Springer - Verlag, 1976.
\item[11a.] Andrew Wiles, Modular elliptic curves and Fermat's Last Theorem, Annals of Mathematics 141 (1995), pp. 443-551.
\item[11b.] R.Taylor and A.Wiles, Ring theoretic properties of certain Hecke algebras, Annals of Mathematics 141 (1995), pp. 553-572.
\end{description}

\end{document}